\documentclass[11pt]{article}

\title{Free field constructions for the elliptic algebra
${\cal A}_{q,p}(\widehat{sl}_2)$ and 
Baxter's eight-vertex model\footnote{
Based on talks given at ``MATHPHYS ODYSSEY 2001--Integrable Models and Beyond''
at Okayama and Kyoto, February 2001, 
EuroConference ``Applications of the
Macdonald Polynomials'' at Newton institute, April 2001, 
ICMS Workshop
``Classical and Quantum Integrable Systems and their Symmetries''
at Heriot-Watt
University, December 2001, 
and  ``Integrable Structures of Exactly Solvable
Two-Dimensional Models of Quantum Field Theory'' at Chernogolovka, 
September 2002.}}

\author{Jun'ichi Shiraishi\\
\\
{\it Graduate School of 
Mathematical Science, }\\ {\it University of Tokyo, Tokyo, Japan}}

\date{}

\begin{document}

\maketitle

\begin{abstract}
Three examples of free field constructions for the vertex operators of
the elliptic quantum group
${\cal A}_{q,p}(\widehat{sl}_2)$ are obtained.
Two of these (for $p^{1/2}=\pm q^{3/2},p^{1/2}=-q^2$) are
based on representation theories of the deformed Virasoro algebra,
which correspond to the level $4$ and level $2$ 
$Z$-algebra of Lepowsky and Wilson.
The third one ($p^{1/2}=q^{3}$) is
constructed over a tensor product of a bosonic and a fermionic
Fock spaces. The algebraic structure at $p^{1/2}=q^{3}$, however,
is not related to the deformed Virasoro algebra.
Using these free field constructions, 
an integral formula for the
correlation functions of
Baxter's eight-vertex model
is obtained. 
This formula shows different structure compared with the
one obtained by Lashkevich and Pugai.

\end{abstract}

%%%%%%%%%%%%%%%%%%%%%%%%%%%%%%%%%%%%%%%%%%%%%%%%%%%%%%%%%%%%%%%%%%%%
%%%%%%%%%%%%%%%%%%%%%%%%%%%%%%%%%%%%%%%%%%%%%%%%%%%%%%%%%%%%%%%%%%%%
\section{Introduction}

Baxter's eight-vertex model is one of the
most fundamental objects
in two dimensional integrable models or in one dimensional integrable
quantum spin chain systems. Basic quantities such as the free energy,
critical exponents,
and excited states were studied almost thirty years ago.
(See \cite{B1} and Refs. in \cite{B}). 
A beautiful formula for the spontaneous polarization 
$P_0$ was obtained by Baxter and Kelland \cite{BK}. It reads
\begin{eqnarray}
P_0 =\prod_{n=1}^\infty
\left(
{1+p^n\over 1-p^n}{1-q^{2n}\over 1+q^{2n}}
\right)^2, \label{spon}
\end{eqnarray}
where $p$ denotes the elliptic nome and $q$ corresponds to the 
crossing parameter for Baxter's elliptic $R$ matrix.

An algebraic approach based on Baxter's corner transfer 
matrix method and the elliptic quantum group for $\widehat{sl}_2$ was
proposed by Kyoto school \cite{JMN}\cite{FIJKMY1},\cite{FIJKMY2}.
They constructed the so-called vertex operator
formalism for the correlation functions for the eight-vertex model
based on the elliptic algebra ${\cal A}_{q,p}(\widehat{sl}_2)$. 
It was shown that 
this elliptic quantum algebra 
${\cal A}_{q,p}(\widehat{sl}_2)$ can be obtained as a quasi-Hopf 
twist of the quantum
affine algebra $U_q(\widehat{sl}_2)$ by a certain twistor element 
\cite{Fr1}\cite{Fr2}\cite{JKOS}.
Lashkevich and Pugai \cite{LaP}\cite{La} succeeded in finding
a bosonic realization of the correlation functions
in terms of the vertex operators 
for the Andrews-Baxter-Forrester (ABF) model \cite{LP}.
This is due to their nice representation of 
a particular kind of the grading operator
for the eight-vertex model in terms of 
the screening current for the ABF model, which enabled them to relate 
these two models.
Recently, Quano improved their method and 
obtained simpler expression for the
correlation function \cite{Q}. 

Lashkevich and Pugai's method is based on the
space of states for the ABF model, and it is not easy, at this moment, to
have a complete construction of the grading operator for the eight-vertex model.
Because of this, 
we still do not have a complete understanding of 
the analytic properties of the vertex
operators. 
When we study 
fusion of the vertex operators, we need
to deal with the analyticity in detail.
To this end, it is desired to have a direct construction of the
vertex operators together with 
the grading operator.

As was pointed out in \cite{coset},
the operator for the elementary particle should satisfy a 
deformed-Virasoro-like commutation relation.
Note, however, this algebra differs from
the usual deformed Virasoro algebra, and
to the author's knowledge, nothing has been obtained 
about this, except for the
double Ising case $p=q^4$.
The free field constructions at $p=q^3$ and $q^6$ presented in this
article, give us another examples.

In Section 2, basic properties for the
vertex operators which we need in this paper are summarized.
In Section 3, integral formulae for 
matrix elements, trace, and correlation functions are given.
In Section 4, free field constructions at $p=q^3,q^4,q^6$ 
are argued.
Most of the statements given in this article have not been proved yet.
Validity of these are examined
by studying matrix elements, trace formulae, or 
$q$-expansions.

\section{Vertex operators for Baxter's eight-vertex model}

%%%%%%%%%%%%%%%%%%%%%%%%%%%%%%%%%%%%%%%%%%%%%%%%%%%%%%%%%%%%%%%%%%%%
\subsection{elliptic $R$ matrix}
Let us fix our notations for Baxter's elliptic $R$ matrix and 
recall basic properties for the
vertex operators for ${\cal A}_{q,p}(\widehat{sl}_2)$ 
\cite{JMN}\cite{FIJKMY1}\cite{FIJKMY2}\cite{JKOS}.
The eight-vertex model is defined by
the elliptic solution 
to the Yang-Baxter equation
\begin{eqnarray}
R_{12}(\zeta_1/\zeta_2)R_{13}(\zeta_1/\zeta_3)R_{23}(\zeta_2/\zeta_3)=
R_{23}(\zeta_2/\zeta_3)R_{13}(\zeta_1/\zeta_3)R_{12}(\zeta_1/\zeta_2),
\end{eqnarray}
obtained by Baxter \cite{B1}. It takes the form as
\begin{eqnarray}
R(\zeta)=
\left(
\begin{array}{cccc}
a(\zeta)&&&d(\zeta)\\
&b(\zeta)&c(\zeta)&\\
&c(\zeta)&b(\zeta)&\\
d(\zeta)&&&a(\zeta)
\end{array}
\right),
\end{eqnarray}
where $a,b,c$ and $d$ are given by elliptic (theta) functions
up to a normalization factor.
For our later purpose, namely for the vertex operator approach,
the $R$ matrix should be normalized in such a way that
the partition function per
site becomes unity. Let us write down $a,b,c$ and $d$
with this normalization. 
We have the parity relations
$a(-\zeta)=-a(\zeta),b(-\zeta)=-b(\zeta),
c(-\zeta)=c(\zeta),d(-\zeta)=d(\zeta)$, and
\begin{eqnarray}
&&
a(\zeta)+d(\zeta)=\zeta^{-1}{\alpha(\zeta^{-1} ) \over \alpha(\zeta)},\\
&&
b(\zeta)+c(\zeta)={\beta(\zeta^{-1}) \over \beta(\zeta)},\\
&&
\alpha(\zeta)
=\sum_{l=0}^\infty \alpha_l \zeta^l=
{(-p^{1/2}q \zeta;p)_\infty \over (-p^{1/2}q^{-1} \zeta;p)_\infty }
{1\over \xi(\zeta^2;p,q)},\\
&&
\beta(\zeta)
=\sum_{l=0}^\infty \beta_l \zeta^l=
{(-q \zeta;p)_\infty \over (-p q^{-1} \zeta;p)_\infty }
{1\over \xi(\zeta^2;p,q)},\\
&&
\xi(z;p,q)=
{(q^2 z;p,q^4)_\infty (pq^2 z;p,q^4)_\infty 
\over (q^4 z;p,q^4)_\infty (p z;p,q^4)_\infty }.
\end{eqnarray}
Here, we used the standard notation 
\begin{eqnarray}
&&
(z;p_1,p_2,\cdots)_\infty=\prod_{i_1,i_2,\cdots=0}^\infty 
(1-p_1^{i_1}p_2^{i_2}\cdots z).
\end{eqnarray}

%%%%%%%%%%%%%%%%%%%%%%%%%%%%%%%%%%%%%%%%%%%%%%%%%%%%%%%%%%%%%%%%%%%%
\subsection{type I vertex operator}

The type I vertex operators (VO's) are defined as intertwiners
between ${\cal A}_{q,p}(\widehat{sl}_2)$-modules
\begin{eqnarray}
\Phi_{\Lambda_i}^{\Lambda_{1-i},V}(\zeta):V(\Lambda_i)\rightarrow 
V(\Lambda_{1-i})\otimes V_\zeta,
\end{eqnarray}
where $V(\Lambda_0),V(\Lambda_1)$ are the elliptic counterpart
of the level one irreducible highest weight modules of
$U_q(\widehat{sl}_2)$, and
$V_\zeta={\bf C}[\zeta,\zeta^{-1}]\otimes({\bf C}v_+\oplus {\bf C}v_-)$ is
the  spin $1/2$ evaluation module.
The components of the vertex operators are defined as
$\Phi_{\Lambda_i}^{\Lambda_{1-i},V}(\zeta)=\Phi_+^{(i)}(\zeta)\otimes v_+ 
+\Phi_-^{(i)}(\zeta)\otimes v_-$.
We adopt the normalization 
$
 \Phi_-^{(0)}(\zeta)|\Lambda_0\rangle=|\Lambda_1\rangle+\cdots
$ and
$
 \Phi_+^{(1)}(\zeta)|\Lambda_1\rangle=|\Lambda_0\rangle+\cdots.
$
These components have the parity relation
$
\Phi_\pm^{(i)}(-\zeta)=\mp (-)^i \Phi_\pm^{(i)}(\zeta).
$

It was shown that the vertex operators enjoy the 
commutation relation
\begin{eqnarray}
\sum_{\epsilon_1',\epsilon_2'=\pm}
R_{\epsilon_1\epsilon_2}^{\epsilon_1'\epsilon_2'}(\zeta_1/\zeta_2)
\Phi_{\epsilon_1'}^{(1-i)}(\zeta_1)\Phi_{\epsilon_2'}^{(i)}(\zeta_2)
=\Phi_{\epsilon_2}^{(1-i)}(\zeta_2)\Phi_{\epsilon_1}^{(i)}(\zeta_1).
\end{eqnarray}
We note, however, that the argument in \cite{JKOS} 
works under the
assumption that the elliptic nome $p$ is
infinitesimally small. So some direct constructions of 
${\cal A}_{q,p}(\widehat{sl}_2)$ (or
$\Phi_{\Lambda_i}^{\Lambda_{1-i},V}(\zeta)$)
for finite $p$ are
desired.

When we try to bosonize the vertex operators, 
we eventually notice that it is convenient to
introduce another basis of
$V_\zeta$
\begin{eqnarray}
&&v_0={1\over 2}(v_+ +v_-),\qquad \qquad v_1={1\over 2}(-v_+ +v_-).
\end{eqnarray}
New components of the vertex operators are 
defined accordingly by
\begin{eqnarray}
&&\Phi_{\Lambda_i}^{\Lambda_{1-i},V}(\zeta)=\Phi_0^{(i)}(\zeta)\otimes v_0
+\Phi_1^{(i)}(\zeta)\otimes v_1,\\
&&
\Phi_0^{(i)}(\zeta)=\Phi_+^{(i)}(\zeta)+\Phi_-^{(i)}(\zeta),\\
&&
\Phi_1^{(i)}(\zeta)=-\Phi_+^{(i)}(\zeta)+\Phi_-^{(i)}(\zeta)=
(-)^i\Phi_0^{(i)}(-\zeta). \label{phi1}
\end{eqnarray}
Since $\Phi_0^{(i)}(\zeta)$ and $\Phi_1^{(i)}(\zeta)$
are related as above, it is enough to consider 
$\Phi_0^{(i)}(\zeta)$ only. 
Our task is to find free field realizations
of $\Phi_0^{(i)}(\zeta)$. 
In what follows, we shall use the shorthand notation
\begin{eqnarray}
\Phi_0^{(i)}(\zeta)\equiv\Phi(\zeta).\label{phi} 
\end{eqnarray}
Then, for example, the original components are written as
\begin{eqnarray}
\Phi_{\pm}^{(i)}(\zeta)={1\over 2}
\left(
\Phi (\zeta)\mp (-)^i \Phi (-\zeta)
\right).
\end{eqnarray}

For calculating matrix elements of the vertex operators,
we need to deal with the Fourier components 
\begin{eqnarray}
\Phi(\zeta)=\sum_{n\in {\bf Z}} \Phi_n \zeta^{-n},
\end{eqnarray}
and commutation relations among them. The correct
relations can be obtained 
after the Riemann-Hilbert splitting of $R\Phi\Phi=\Phi\Phi$ relation.
 These are
\begin{eqnarray}
\zeta_1^{-1}\alpha(\zeta_2/\zeta_1)\left(
\Phi(\zeta_1)\Phi(\zeta_2)-\Phi(-\zeta_1)\Phi(-\zeta_2)
\right)
&=&(1\leftrightarrow 2),\\
\beta(\zeta_2/\zeta_1)\left(
\Phi(\zeta_1)\Phi(\zeta_2)+\Phi(-\zeta_1)\Phi(-\zeta_2)
\right)
&=&(1\leftrightarrow 2).
\end{eqnarray}
If we write them in modes, we have
\begin{eqnarray}
&&\sum_{l=0}^\infty \alpha_l \Phi_{n-l}\Phi_{m+l}=
\sum_{l=0}^\infty \alpha_l \Phi_{m-l-1}\Phi_{n+l+1}\qquad\qquad
(n+m:{\rm odd}),\\
&&\sum_{l=0}^\infty \beta_l \Phi_{n-l}\Phi_{m+l}=
\sum_{l=0}^\infty \beta_l \Phi_{m-l}\Phi_{n+l}\qquad\qquad\qquad
(n+m:{\rm even}).
\end{eqnarray}

Using the commutation relations and the `highest weight conditions'
\begin{eqnarray}
&&\Phi_0|\Lambda_i\rangle=|\Lambda_i\rangle,\qquad
\langle \Lambda_i|\Phi_{0}=\langle \Lambda_i|,
\qquad \langle \Lambda_i|\Lambda_i\rangle=1,\\
&&\Phi_n|\Lambda_i\rangle=0,\qquad\langle \Lambda_i|\Phi_{-n}=0
\qquad(n>0), \nonumber
\end{eqnarray}
we are able to calculate any matrix
elements, in principle. For example, 
we have
\begin{eqnarray}
&&
\langle \Lambda_i| \Phi(\zeta_1)\Phi(\zeta_2)|\Lambda_i\rangle=
\beta(\zeta_2/\zeta_1)^{-1},\\
&&
\langle \Lambda_{1-i}| \Phi_1\Phi(\zeta_1)\Phi(\zeta_2)|\Lambda_i \rangle=
-\beta_1 \zeta_1 \alpha(\zeta_2/\zeta_1)^{-1},
\end{eqnarray}
and so on.

If we define dual of 
the components by
\begin{eqnarray}
&&
\Phi^{*(i)}_0(\zeta)={1\over 2}\xi(1;p,q) \Phi^{(i)}_0(-q^{-1}\zeta), 
\label{phi*1}\\
&&
\Phi^{*(i)}_1(\zeta)=-{1\over 2}\xi(1;p,q)
\Phi^{(i)}_1(-q^{-1}\zeta),\label{phi*2}
\end{eqnarray}
we have the inversion relation
\begin{eqnarray}
 \Phi^{*(1-i)}_0(\zeta)\Phi^{(i)}_0(\zeta)
+
 \Phi^{*(1-i)}_1(\zeta)\Phi^{(i)}_1(\zeta)
={\rm id}.
\end{eqnarray}

%%%%%%%%%%%%%%%%%%%%%%%%%%%%%%%%%%%%%%%%%%%%%%%%%%%%%%%%%%%%%%%%%%%%
\subsection{type II vertex operator}
The type II vertex operators satisfy
\begin{eqnarray}
-\Psi_{\epsilon_1}^{(1-i)}(\zeta_1)\Psi_{\epsilon_2}^{(i)}(\zeta_2)
=
\sum_{\epsilon_1',\epsilon_2'=\pm}
\Psi_{\epsilon_2'}^{(1-i)}(\zeta_2)\Psi_{\epsilon_1'}^{(i)}(\zeta_1)
{R^*}_{\epsilon_1\epsilon_2}^{\epsilon_1'\epsilon_2'}(\zeta_1/\zeta_2),
\end{eqnarray}
where the matrix $R^*(\zeta)$ is given 
by shifting the elliptic nome in the $R$ matrix as 
\begin{eqnarray}
R^*(\zeta)=R(\zeta)\Bigl|_{p^{1/2}\rightarrow p^{*1/2}}.
\end{eqnarray}
where $p^{*1/2}=p^{1/2}q^{-1}$.
Note that at the series of points
\begin{eqnarray}
p^{1/2}=\pm q^{1+1/l}\qquad \qquad (l=1,2,3,\cdots),
\end{eqnarray}
the $R^*$ becomes 
(anti-)diagonal matrix. Correspondingly,
the structure of the algebra for the type II VO
becomes simpler.
For $p=q^4$, $\Psi(\zeta)$ reduces to a free fermion.
For $p=q^3$,
$\Psi(\zeta)$ is no more a free field, nevertheless it
can be described by the deformed Virasoro algebra
with a special choice of the parameters. (See Section \ref{freeboson}.)

Introducing the sum of the components
\begin{eqnarray}
\Psi(\zeta)=\Psi^{(i)}_+(\zeta)+\Psi^{(i)}_-(\zeta),
\end{eqnarray}
as we did for the type I VO, we have the 
exchange relation for $\Psi(\zeta)$ as follows: 
\begin{eqnarray}
&&\zeta_1^{-1}{1\over 1-\zeta_2^2/\zeta_1^2}
{1\over \alpha^*(\zeta_2/\zeta_1)}\left(
\Psi(\zeta_1)\Psi(\zeta_2)-\Psi(-\zeta_1)\Psi(-\zeta_2)
\right) \nonumber\\
&&=(1\leftrightarrow 2),\\
&&{(1+q^{-1}\zeta_2/\zeta_1)(1+q \zeta_2/\zeta_1)\over
 1-\zeta_2^2/\zeta_1^2}
{1\over \beta^*(\zeta_2/\zeta_1)}\left(
\Psi(\zeta_1)\Psi(\zeta_2)+\Psi(-\zeta_1)\Psi(-\zeta_2)
\right) \nonumber\\
&&=(1\leftrightarrow 2),
\end{eqnarray}
where 
$\alpha^*(\zeta)=\alpha(\zeta)\Bigl|_{p^{1/2}\rightarrow p^{*1/2}}$
and 
$\beta^*(\zeta)=\beta(\zeta)\Bigl|_{p^{1/2}\rightarrow p^{*1/2}}$.

The relations between the type I and type II VO is
the scalar exchange relation
\begin{eqnarray}
\zeta_1{(q\zeta_2^2/\zeta_1^2;q^4)_\infty \over 
(q^3\zeta_2^2/\zeta_1^2;q^4)_\infty }
\Phi(\zeta_1)\Psi(\zeta_2)=
\zeta_2{(q\zeta_1^2/\zeta_2^2;q^4)_\infty \over 
(q^3\zeta_1^2/\zeta_2^2;q^4)_\infty
}
\Psi(\zeta_2)\Phi(\zeta_1). \label{phipsi}
\end{eqnarray}
At $p=q^3,q^4$, $\Psi(\zeta)$ can be realized by 
the deformed Virasoro current $T(\zeta)$.
When we investigate free field formulae for $\Phi(\zeta)$,
the scalar exchange relation (\ref{phipsi}) helps us.

%%%%%%%%%%%%%%%%%%%%%%%%%%%%%%
\subsection{elementary scalar particle for the eight-vertex model}

In the paper \cite{coset}, it was conjectured that 
the operator
\begin{eqnarray}
t(\zeta)&=&C_{II}{\rm Res
}_{\zeta_1=-p^{*1/2}q^{-1}\zeta_2}\Psi(\zeta_1)\Psi(\zeta_2)
{d\zeta_1\over \zeta_1}\label{t}\\
&=&
C_I 
\Phi(\zeta_1)\Phi(\zeta_2)\Bigl|_{\zeta_1=-p^{1/2}q\zeta_2}, \nonumber\\
t(-\zeta)&=&-t(\zeta),
\end{eqnarray}
satisfies a similar commutation relation 
to the deformed Virasoro algebra, namely
\begin{eqnarray}
&&f\left({\zeta_2^2\over \zeta_1^2}\right)t(\zeta_1)t(\zeta_2)=
f\left({\zeta_1^2\over  \zeta_2^2}\right)t(\zeta_2)t(\zeta_1)+ \label{comfort}\\[3mm]
&&+
{(1-p)(1-p^{*-1})\over 1-q^2}{1\over 2}
\left(
\delta\left({q\zeta_2\over \zeta_1}\right)-
\delta\left({\zeta_2\over q\zeta_1}\right)-
\delta\left(-{q\zeta_2\over \zeta_1}\right)+
\delta\left(-{\zeta_2\over q\zeta_1}\right)
\right),\nonumber
\end{eqnarray}
where
\begin{eqnarray}
f(z)={1\over 1-z}
{(pz;q^4)_\infty\over (pq^2z;q^4)_\infty}
{(p^{*-1}z;q^4)_\infty\over (p^{*-1}q^2z;q^4)_\infty},
\end{eqnarray}
and $\delta(\zeta)=\sum_{n\in{\bf Z}} \zeta^n$.
This operator $t(\zeta)$ 
represent the Fadeev-Zamolodchikov
algebra for the elementary scalar particles 
for the eight-vertex model.

Contrary to the case of the
ordinary deformed Virasoro algebra,
to obtain an explicit realization of this 
algebra seems a not easy task. 
In \cite{coset}, $t(\zeta)$ for 
the double Ising case ($p=q^4$) was presented in terms
of a free fermion.
Up to now,
however, 
no other example of free field construction has been obtained,
to the author's knowledge.
We can check that the free field constructions for
$\Phi(\zeta)$ (or $\Psi(\zeta)$) at $p=q^3,q^4,q^6$
given in Section 4 provide us with 
other examples of explicit
formulae for $t(\zeta)$.

%%%%%%%%%%%%%%%%%%%%%%%%%%%%%%%%%%%%%%%%%%%%%%%%%%%%%%%%%%%%%%%%%%%%
%%%%%%%%%%%%%%%%%%%%%%%%%%%%%%%%%%%%%%%%%%%%%%%%%%%%%%%%%%%%%%%%%%%%
\section{Results}
Integral formulae for matrix elements,
trace over the irreducible highest weight modules,
and correlation functions for the inhomogeneous eight-vertex model
are given in this section.
These are obtained from the use of free field constructions
given in the next section.

%%%%%%%%%%%%%%%%%%%%%%%%%%%%%%%%%%%%%%%%%%%%%%%%%%%%%%%%%%%%%%%%%%%%
\subsection{matrix elements}

By using the free field constructions of 
the type I vertex operator 
$\Phi(\zeta)=\Phi_0^{(i)}(\zeta)$ (see Eq. (\ref{phi}))
for $p^{1/2}=q^{3/2},-q^2,q^3$, we have a
integral formula for the matrix elements.
Introduce functions as
\begin{eqnarray}
&&
h(\zeta)={(p^{1/2}q^{-1}\zeta;p^{1/2})_\infty\over
(q \zeta;p^{1/2})_\infty}\xi(\zeta^2;p,q),\\
&&
g(\zeta)={(p^{1/4}q^{1/2}\zeta;p^{1/2})_\infty\over
(p^{1/4}q^{-1/2}\zeta;p^{1/2})_\infty}.
\end{eqnarray}
For $n=2,4,6,\cdots$, we have
\begin{eqnarray}
&&\langle \Lambda_i|\Phi(\zeta_1)\Phi(\zeta_2)\cdots 
\Phi(\zeta_n)|\Lambda_i\rangle  \nonumber\\
&=&\left( {(p^{1/2}q^{-1};p^{1/2})_\infty \over (-p^{1/2}q^{-1};p^{1/2})_\infty }
{(p^{1/2} ;p^{1/2})_\infty \over (-p^{1/2};p^{1/2})_\infty }
\right)^n  \prod_{i<j} h(\zeta_j/\zeta_i) \label{matel}\\
&\times& 
\oint\cdots\oint {d\xi_1\over 2 \pi i\xi_1}\cdots {d\xi_n\over 2 \pi i\xi_n}
\prod_{i=1}^n 
{\Theta_{p^{1/2}}(-p^{1/4}q^{-1/2}\zeta_i/\xi_i) \over
\Theta_{p^{1/2}}(p^{1/4}q^{-1/2}\zeta_i/\xi_i) } \nonumber\\
&\times&
\prod_{k=1}^n 
\left[
\prod_{i<k}g(\zeta_k/\xi_i)\prod_{j\geq k}g(\xi_j/\zeta_k)
\right] 
F(\xi_1,\cdots,\xi_n;p^{1/2},q), \nonumber
\end{eqnarray}
where the integration contour for $\xi_i$ is 
given by the condition $|\xi_i/\zeta_i|=1$, namely,
that encloses
poles at $\xi_i=p^{1/4+m/2}q^{-1/2}\zeta_i$ ($m=0,1,2\cdots$), and 
the function $F(\xi_1,\cdots,\xi_n;p^{1/2},q)$ is defined by
\begin{eqnarray}
&&
F(\xi_1,\cdots,\xi_n;q^{3/2},q)
=\prod_{i<j}
(1-\xi_j/\xi_i)
{(-q\xi_j/\xi_i;q^{3/2})_\infty \over (-q^{1/2}\xi_j/\xi_i;q^{3/2})_\infty },\\
&&
F(\xi_1,\cdots,\xi_n;-q^{2},q)
=\prod_{i<j \atop {\rm step\;2}}
(1-\xi_j/\xi_i)
{(-q\xi_j/\xi_i;-q^2)_\infty \over (q \xi_j/\xi_i;-q^2)_\infty },\\
&&
F(\xi_1,\cdots,\xi_n;q^{3},q)
=
{\rm Pfaffian}\left({1-\xi_j^2/\xi_i^2\over 
(1+q^{-1}\xi_j/\xi_i)(1+q\xi_j/\xi_i)} \right)_{1\leq i,j\leq n}\nonumber\\
&&\qquad\qquad\qquad \qquad\qquad \qquad\times
\prod_{i<j} 
{(-q\xi_j/\xi_i;q^{3})_\infty \over (-q^{2}\xi_j/\xi_i;q^{3})_\infty },
\end{eqnarray}
where we used the notation
$\prod_{i<j \atop {\rm step\;2}} f_{ij}=f_{13}f_{15}\cdots
f_{24}f_{26}\cdots$.

It is amusing to see that the matrix elements
can be written down in a similar manner for different
$p^{1/2}$.
At present, we have different constructions
for $\Phi(\zeta)$ for each $p=q^3$, $p=q^4$ and $p=q^6$,
and there is no a priori reason to have 
similar integral formulae.
To investigate the reason for this phenomenon
might help our further study on the vertex operators for
general $p^{1/2}$.

To give an heuristic argument for how one can
obtain the formula (\ref{matel}) is in order.
This formula was first obtained at $p=q^3$ and $p=q^4$,
where we have constructions for 
$\Phi(\zeta)$ based on the representation theory of the 
deformed Virasoro algebra.
As it will be explained in the next section, the
bosonic constructions for $p=q^3$ and $p=q^4$ looks quite
different.
One realizes, however, that the resulting integral formulae for 
the matrix elements share the same structures for
the part given by the functions
$h(\zeta)$ and $g(\zeta)$. 
Moreover, for the two point case ($n=2$), the function
$F(\xi_1,\xi_2;p^{1/2},q)$ can be extrapolated for general $p^{1/2}$
(see Eq. (\ref{Ffor2})).
Then we realize a nice simplification of
$F(\xi_1,\xi_2;p^{1/2},q)$ as $p^{1/2}=q^3$ (see Eq.(\ref{Fq6})).
This gives us be a hint for a free field construction for 
$\Phi$ at $p^{1/2}=q^3$. Thus we arrive at (\ref{matel}).

%%%%%%%%%%%%%%%%%%%%%%%%%%%%%%%%%%%%%%%%%%%%%%%%%%%%%%%%%%%%%%%%%%%%
\subsection{two point matrix elements for 
general $p^{1/2}$}

For the two point matrix element
\begin{eqnarray}
\langle \Lambda_i|\Phi(\zeta_1)\Phi(\zeta_2)|\Lambda_i\rangle=
{1\over \beta(\zeta_2/\zeta_1)},\label{2pt}
\end{eqnarray}
the integration kernel $F(\xi_1,\xi_2;p^{1/2},q)$ 
for general $p^{1/2}$ 
can be investigated and written down in a simple form:
\begin{eqnarray}
&&F(\xi_1,\xi_2;p^{1/2},q) \nonumber\\
&=&
(1-\xi_2/\xi_1)\; {}_2\phi_1\left({p^{1/2}q^{-1},-pq^{-2}\atop-q};
p^{1/2},-p^{-1/2}q^2 \xi_2/\xi_1\right) \label{Ffor2}\\
&=&\sum_{n=0}^\infty c_n \; (\xi_2/\xi_1)^n.\nonumber
\end{eqnarray}
Here ${~}_2\phi_1\left({a,b\atop c};q,z\right)$ denotes the 
basic hypergeometric series defined by
\begin{eqnarray}
&&{~}_2\phi_1\left({a,b\atop c};q,z\right)
=\sum_{n=0}^\infty 
{(a;q)_n(b;q)_n\over (c;q)_n(q;q)_n}z^n,
\end{eqnarray}
where $(a;q)_n=(1-a)(1-aq)\cdots(1-aq^{n-1})$.

The coefficients $c_n$ can be written as
\begin{eqnarray}
&&c_n\\
&=&{(q^{-1};p^{1/2})_n(p^{1/2}q^{-1/2};p^{1/2})_n
(-p^{1/2}q^{-1/2};p^{1/2})_n(-p^{1/2}q^{-2};p^{1/2})_n
\over 
(-q;p^{1/2})_n( q^{-1/2};p^{1/2})_n
(- q^{-1/2};p^{1/2})_n(p^{1/2} ;p^{1/2})_n}
\;(-p^{-1/2}q^2)^n. \nonumber
\end{eqnarray}
After integrating (\ref{matel}) for $n=2$ with (\ref{Ffor2}), we have
the identity
\begin{eqnarray}
&&\sum_{k=0}^\infty
{}_2\phi_1\left({q,-qp^{k/2}\atop -p^{(k+1)/2}};
p^{1/2},p^{1/2}q^{-1}\zeta\right)
\;{}_2\phi_1\left({q,-p^{k/2}\atop -p^{(k+1)/2}q^{-1}};
p^{1/2},p^{1/2}q^{-1}\zeta\right)\times \nonumber\\
&&
\qquad\qquad
\times {(-q;p^{1/2})_k(-1;p^{1/2})_k\over 
(-p^{1/2};p^{1/2})_k (-p^{1/2}q^{-1};p^{1/2})_k}p^{k/2}q^{-k}c_k\\
&=&
{(-pq^{-1}\zeta;p)_\infty \over (-q \zeta;p)_\infty}
{(q\zeta;p^{1/2})_\infty \over (p^{1/2}q^{-1} \zeta;p^{1/2})_\infty}. \nonumber
\end{eqnarray}
This means (\ref{2pt}).

At $p^{1/2}=q^3$,
the function $F(\xi_1,\xi_2;p^{1/2},q)$ 
shows special degeneration, namely we have 
\begin{eqnarray}
F(\xi_1,\xi_2;q^3,q)=(1-\xi_2^2/\xi_1^2)
{(-q^4\xi_2/\xi_1;q^3)_\infty \over(-q^{-1}\xi_2/\xi_1;q^3)_\infty }.
\label{Fq6}
\end{eqnarray}
This helps us when we study free field realization of 
$\Phi(\zeta)$ at $p^{1/2}=q^3$.

%%%%%%%%%%%%%%%%%%%%%%%%%%%%%%%%%%%%%%%%%%%%%%%%%%%%%%
\subsection{trace over irreducible representations}

Since we have free field formulae for $\Phi(\zeta)$ 
at $p=q^3,q^4,q^6$,
to calculate 
trace over irreducible highest weight modules is 
a straightforward task. 

Let $x$ be a parameter satisfying $|x|<1$.
Introduce
\begin{eqnarray}
&&
\tilde{h}(\zeta)=\prod_{k=0}^\infty h(x^k\zeta)
\prod_{k=1}^\infty h(x^k\zeta^{-1}),\\
&&
\tilde{g}(\zeta)=\prod_{k=0}^\infty g(x^k\zeta)
\prod_{k=1}^\infty g(x^k\zeta^{-1}),
\end{eqnarray}
and for $n=2,4,6,\cdots$, we define
\begin{eqnarray}
&&
\tilde{F}(\xi_1,\cdots,\xi_n;q^{3/2},q,x)
=\prod_{i<j}
(\xi_j/\xi_i;x)_\infty(x\xi_i/\xi_j;x)_\infty\times\\
&&\qquad\qquad\qquad\times
{(-q\xi_j/\xi_i;q^{3/2},x)_\infty \over
(-q^{1/2}\xi_j/\xi_i;q^{3/2},x)_\infty }
{(-qx\xi_i/\xi_j;q^{3/2},x)_\infty \over
(-q^{1/2}x\xi_i/\xi_j;q^{3/2},x)_\infty },\nonumber\\ 
&&
\tilde{F}(\xi_1,\cdots,\xi_n;-q^{2},q,x)
=\prod_{i<j \atop {\rm step\;2}}
(\xi_j/\xi_i;x)_\infty(x\xi_i/\xi_j;x)_\infty\times\\
&&\qquad\qquad\qquad\times
{(-q\xi_j/\xi_i;-q^2,x)_\infty \over (q \xi_j/\xi_i;-q^2,x)_\infty }
{(-qx\xi_i/\xi_j;-q^2,x)_\infty \over (qx \xi_i/\xi_j;-q^2,x)_\infty },
\nonumber\\
&&
\tilde{F}(\xi_1,\cdots,\xi_n;q^{3},q,x)\label{Fforq^6}\\
&&\qquad
=
{\rm Pfaffian}\left(
{ (x;x)_\infty^2\over 2(-x;x)_\infty^2}
\left(
{\Theta_x(q\xi_j/\xi_i)\over\Theta_x(-q\xi_j/\xi_i)}+
{\Theta_x(q^{-1}\xi_j/\xi_i)\over\Theta_x(-q^{-1}\xi_j/\xi_i)}
\right)
\right)_{1\leq i,j\leq n} \nonumber\\
&&\qquad\qquad\qquad\times
\prod_{i<j} 
{(-q\xi_j/\xi_i;q^3,x)_\infty \over (-q^2\xi_j/\xi_i;q^3,x)_\infty}
{(-qx\xi_i/\xi_j;q^3,x)_\infty \over (-q^2x\xi_i/\xi_j;q^3,x)_\infty},
\nonumber\\
&&
E(\zeta_1,\cdots,\zeta_n,\xi_1,\cdots,\xi_n;q^{3/2},q,x)
=
{\Theta_{x^{6}} \left(x^4\prod_{i=1}^n {\zeta_i^2\over \xi_i^2}\right)
\over 
\Theta_{x^6}(x^4)} ,\\
&&
E(\zeta_1,\cdots,\zeta_n,\xi_1,\cdots,\xi_n;-q^{2},q,x)\\
&&\qquad \qquad 
=
{\Theta_{x^{4}} \left(x^2\prod_{i=1}^{n/2} 
{\zeta_{2i-1}\zeta_{2i}\over\xi_{2i-1}^2}\right)
\Theta_{x^{4}} \left(x^3\prod_{i=1}^{n/2} 
{\zeta_{2i-1}\zeta_{2i}\over\xi_{2i}^2}\right)\over
\Theta_{x^4}(x^2)\Theta_{x^4}(x^3)} ,\nonumber\\ 
&&
E(\zeta_1,\cdots,\zeta_n,\xi_1,\cdots,\xi_n;q^{3},q,x)
=
{\Theta_{x^{3}} \left(x^2\prod_{i=1}^n {\zeta_i\over \xi_i}\right) 
\over \Theta_{x^3}(x^2)}.
\end{eqnarray}
Set further
\begin{eqnarray}
&&C_\phi=\prod_{k=1}^\infty h(x^k),\\
&&C_S=
\left\{
\begin{array}{ll}
(x;x)_\infty{\displaystyle (-qx;q^{3/2},x)_\infty \over
\displaystyle (-q^{1/2}x;q^{3/2},x)_\infty } \qquad
& (p^{1/2}=q^{3/2}),\\[4mm]
(x;x)_\infty{\displaystyle (-qx;-q^2,x)_\infty \over
\displaystyle (q x;-q^2,x)_\infty } \qquad
& (p^{1/2}=-q^2),\\[4mm]
 {\displaystyle (-qx;q^3,x)_\infty \over
\displaystyle (-q^2 x;q^3,x)_\infty } \qquad
& (p^{1/2}=q^3),
\end{array}\right..
\end{eqnarray}

We have the 
integral representation for the trace  which works for
$p^{1/2}=q^{3/2},-q^2,q^3$:
\begin{eqnarray}
&&
{1\over 
{\rm tr}_{V(\Lambda_i)}(x^D)}
{\rm tr}_{V(\Lambda_i)}
\Bigl(x^{D}\Phi(\zeta_1)\Phi(\zeta_2)\cdots 
\Phi(\zeta_n)\Bigr)\nonumber\\
&=& 
\left( C_\phi C_S{(p^{1/2}q^{-1};p^{1/2})_\infty
\over (-p^{1/2}q^{-1};p^{1/2})_\infty } {(p^{1/2} ;p^{1/2})_\infty \over
(-p^{1/2};p^{1/2})_\infty }
\right)^n  \prod_{i<j} \tilde{h}(\zeta_j/\zeta_i)\times \nonumber\\
&\times& 
\oint\cdots\oint {d\xi_1\over 2 \pi i\xi_1}\cdots {d\xi_n\over 2 \pi i\xi_n}
\prod_{i=1}^n 
{\Theta_{p^{1/2}}(-p^{1/4}q^{-1/2}\zeta_i/\xi_i) \over
\Theta_{p^{1/2}}(p^{1/4}q^{-1/2}\zeta_i/\xi_i) } \times \label{trace} \\
&\times&
\prod_{k=1}^n 
\left[
\prod_{i<k}\tilde{g}(\zeta_k/\xi_i)\prod_{j\geq k}\tilde{g}(\xi_j/\zeta_k)
\right] 
\tilde{F}(\xi_1,\cdots,\xi_n;p^{1/2},q,x) \times\nonumber\\
&\times&
E(\zeta_1,\cdots,\zeta_n,\xi_1,\cdots,\xi_n;p^{1/2},q,x),\nonumber
\end{eqnarray}
where $n=2,4,6,\cdots$, $D$ is the corner Hamiltonian, and integration contour
for $\xi_i$ encloses poles at $\xi_i=x^l p^{1/4+m/2}q^{-1/2}\zeta_i$
($l=1,2,\cdots,m=0,1,2\cdots$).

%%%%%%%%%%%%%%%%%%%%%%%%%%%%%%%%%%%%%%%%%%%%%%%%%%%%%%%
\subsection{correlation functions for the eight-vertex model}
In view of our free field formula,
we need to represent local operators for the
basis $v_0,v_1$.
For example, $\sigma_z$ (which satisfies $\sigma_z v_\pm=\pm v_\pm$) is represented as
\begin{eqnarray}
\sigma_z=\left(\begin{array}{cc} 0&-1\\-1&0 \end{array}\right).
\end{eqnarray}
By $E_{\eta' \eta}$ ($\eta,\eta'=0,1$) we denote
the matrix unit with respect to the basis $v_0,v_1$.
Using (\ref{phi1}), (\ref{phi}), (\ref{phi*1}), (\ref{phi*2}), 
and setting $x=q^2$ in  (\ref{trace}), we
have the following integral representation of
the correlation function 
\begin{eqnarray}
&&\langle 
E_{\eta_n'\eta_n}\otimes\cdots \otimes E_{\eta_1'\eta_1}
\rangle_{i}\nonumber\\
&=&
{1\over {\rm tr}_{V(\Lambda_i)} (q^{2D})}
 {\rm tr}_{V(\Lambda_i)}
\left(
q^{2D} \Phi^{*(i+1)}_{\eta_1'}(\zeta_1)
\Phi^{*(i+2)}_{\eta_2'}(\zeta_2)\cdots
\Phi^{*(i+n)}_{\eta_n'}(\zeta_n)\times \right.\nonumber\\
&&\qquad\qquad \left. \times
\Phi^{(i+n-1)}_{\eta_n}(\zeta_n)
\cdots
\Phi^{(i+1)}_{\eta_2}(\zeta_2)
\Phi^{(i)}_{\eta_1}(\zeta_1)
\right)\nonumber\\
&=&
(-1)^{\sum_{k=1}^n(k-1+i)(\eta_k'+\eta_k)}
{1\over 2^n}
\left({(q^2;q^2)_\infty\over (p;p)_\infty}\right)^n\times\nonumber\\
& \times&
\left( C_S{(q;p^{1/2})_\infty \over (q;q^2)_\infty}
{(p^{1/2}q^{-1};p^{1/2})_\infty
\over (-p^{1/2}q^{-1};p^{1/2})_\infty } {(p^{1/2} ;p^{1/2})_\infty \over
(-p^{1/2};p^{1/2})_\infty }
\right)^{2n}\times\nonumber\\
&\times&
\prod_{k=1}^n 
\left(
{(-(-)^{\eta_k'+\eta_k}p^{1/2};p^{1/2})_\infty
\over 
(-(-)^{\eta_k'+\eta_k}q^2;q^2)_\infty
}
\right)^2\times\label{correlation}\\
&\times &
\prod_{i<j} 
\tilde{h}\left((-)^{\eta_i'+\eta_j'}{\zeta_j\over\zeta_i}\right)
\tilde{h}\left(-q(-)^{\eta_i'+\eta_j}{\zeta_j\over\zeta_i}\right)
\tilde{h}\left(-q(-)^{\eta_i+\eta_j'}{\zeta_i\over\zeta_j}\right)
\tilde{h}\left((-)^{\eta_i+\eta_j}{\zeta_i\over\zeta_j}\right)\times
\nonumber\\
&\times& \oint\cdots\oint
{d\xi_1'\over 2 \pi i\xi_1'}\cdots {d\xi_n'\over 2 \pi i\xi_n'}
 {d\xi_1\over 2 \pi i\xi_1}\cdots {d\xi_n\over 2 \pi i\xi_n}
\times\nonumber\\
&\times&
\prod_{i=1}^n 
{\Theta_{p^{1/2}}(p^{1/4}q^{-1/2}(-)^{\eta_i'}\zeta_i/\xi_i') \over
\Theta_{p^{1/2}}(-p^{1/4}q^{-1/2}(-)^{\eta_i'}\zeta_i/\xi_i') }
{\Theta_{p^{1/2}}(-p^{1/4}q^{-1/2}(-)^{\eta_i}\zeta_i/\xi_i) \over
\Theta_{p^{1/2}}(p^{1/4}q^{-1/2}(-)^{\eta_i}\zeta_i/\xi_i) } \times \nonumber\\
&\times&
\prod_{k=1}^n 
\left[
\prod_{i<k}
\tilde{g}\left(-(-)^{\eta_k'}{\zeta_k\over\xi_i'}\right)
\prod_{j\geq k}
\tilde{g}\left(-(-)^{\eta_k'}{\xi_j'\over\zeta_k}\right)
\right] \times\nonumber\\
&\times&
\prod_{k=1}^n 
\left[
\prod_{i>k}
\tilde{g}\left((-)^{\eta_k}{\zeta_k\over\xi_i}\right)
\prod_{j\leq k}
\tilde{g}\left((-)^{\eta_k}{\xi_j\over\zeta_k}\right)
\right] \times\nonumber\\
&\times&
\prod_{k=1}^n \prod_{j=1}^n
\left[
\tilde{g}\left(q(-)^{\eta_k}{\zeta_k\over\xi_i'}\right)
\tilde{g}\left(-q(-)^{\eta_k'}{\xi_j\over\zeta_k}\right)
\right] \times\nonumber\\
&\times&
\tilde{F}(q^{-1}\xi_1',\cdots,q^{-1}\xi_n',
\xi_n,\cdots,\xi_1;p^{1/2},q,q^2) \times\nonumber\\
&\times&
E(-q^{-1}(-)^{\eta_1'}\zeta_1,\cdots,-q^{-1}(-)^{\eta_n'}\zeta_n,
(-)^{\eta_n}\zeta_n,\cdots,(-)^{\eta_1}\zeta_1,\nonumber\\
&&\qquad\qquad\qquad\qquad\qquad\qquad\qquad
q^{-1}\xi_1',\cdots,q^{-1}\xi_n'\xi_n,\cdots,\xi_1;p^{1/2},q,q^2),\nonumber
\end{eqnarray}
where $p^{1/2}=q^{3/2},-q^2,q^3$, $1<|\zeta_i|<|q^{-1}|$ and
the integration contours are given by the conditions
$|\zeta_i/\xi_i|=1,|\zeta_i/\xi_i'|=1$.

Notice that for the case $p^{1/2}=q^3$
the factor in the Pfaffian (see (\ref{Fforq^6}))
reduces to the delta function as
\begin{eqnarray}
{(x;x)_\infty^2\over 2(-x;x)_\infty^2}
\left(
{\Theta_x(q\xi_j/\xi_i)\over\Theta_x(-q\xi_j/\xi_i)}+
{\Theta_x(q^{-1}\xi_j/\xi_i)\over\Theta_x(-q^{-1}\xi_j/\xi_i)}
\right)
\Biggl|_{x=q^2}=
\delta\left(
-q^{-1}{\xi_j\over \xi_i}
\right).
\end{eqnarray}
Because of this, half of the integrals, i.e. $n$ integrals 
can be performed.

For the simplest case $n=1$, we can check 
(by $q$-expansion) that
Eq. (\ref{correlation}) correctly gives 
\begin{eqnarray}
&&\Bigl\langle( E_{00}+E_{11}) \Bigr\rangle_i=1,\\
&&\Bigl\langle( -E_{01}-E_{10} )\Bigr\rangle_i=-(-)^i P_0,
\end{eqnarray}
where $P_0$ denotes the Baxter-Kelland formula for the
spontaneous polarization
(\ref{spon}).

%%%%%%%%%%%%%%%%%%%%%%%%%%%%%%%%%%%%%%%%%%%%%%%%%%%%%%
%%%%%%%%%%%%%%%%%%%%%%%%%%%%%%%%%%%%%%%%%%%%%%%%%%%%%%
\section{Free field realizations for the vertex operators}
\label{freeboson}

We list the free field formulae for the vertex operators
$\Phi(\zeta)$ and $\Psi(\zeta)$ for the eight-vertex model.
For $p^{1/2}=\pm q^{3/2}$ and  $p^{1/2}=-q^{2}$
we need to recall the
representation theory of the deformed Virasoro algebra
\cite{SKAO,HJKOS}.
For $p^{1/2}=q^{3}$, a free bosonic and 
a free fermionic field are introduced to realize the VO.

%%%%%%%%%%%%%%%%%%%%%%%%%%%%%%%%%%%%%%%%%%%%%%%%%%%%%%
\subsection{CFT limit}
Let us start from a simple case.
If we set $p^{1/2}=q^r$ and take the limit $q\rightarrow 1$
($r$ is fixed), the functions $\alpha(\zeta)$ and $\beta(\zeta)$
reduce to
\begin{eqnarray}
&&
\alpha(\zeta)\rightarrow {1\over 1+\zeta}
\left(1-\zeta\over 1+\zeta
\right)^{-{r-1\over 2r}},\qquad
\beta(\zeta)\rightarrow 
\left(1-\zeta\over 1+\zeta
\right)^{-{r-1\over 2r}}.
\end{eqnarray}
In this limit, we have a simple bosonic realization
of $\Phi(\zeta)$ and $\Psi(\zeta)$. 
Let $a_n$ be bosons with odd modes satisfying
\begin{eqnarray}
[a_n,a_m]=n\delta_{n+m,0}\qquad  (n,m:{\rm odd}).
\end{eqnarray}
Using this `two-reduced boson', we have
\begin{eqnarray}
&&\Phi(\zeta)=:\exp\left(\sum_{n:{\rm odd}}\sqrt{r-1\over r} {a_n\over n}
\zeta^{-n}\right):, \\
&&
\Psi(\zeta)=:\exp\left(-\sum_{n:{\rm odd}}\sqrt{r\over r-1} {a_n\over n}
\zeta^{-n}\right):.
\end{eqnarray}
One can check that the generators
\begin{eqnarray}
L_n={1\over 4}:\sum_{l:{\rm odd}}a_{2n-l}a_l:+{1\over 16}\delta_{n,0}
\qquad (n\in {\bf Z}),
\end{eqnarray}
satisfy $c=1$ Virasoro algebra, and that
$\zeta^{-(r-1)/2r}\Phi(\zeta)$,
$\zeta^{-r/2(r-1)}\Psi(\zeta)$ are primary fields
to this $c=1$ Virasoro algebra.

At this CFT limit, however, things are oversimplified and 
this bosonic formula
does not help us very much
when we investigate the situation $|q|<1$.
As a matter of fact, it is necessary 
to work on much bigger Fock space and 
reduce the space of states by introducing BRST-type cohomology.

%%%%%%%%%%%%%%%%%%%%%%%%%%%%%%%%%%%%%%%%%%%%%%%%%%%%%%
\subsection{deformed Virasoro algebra as a quantum analogue of Lepowsky-Wilson's
$Z$-algebra}

It was discussed in \cite{HJKOS} that the deformed Virasoro algebra
can be regarded as a smooth deformation 
of Lepowsky-Wilson's $Z$-algebra.
For level $k$, the Virasoro central charge for 
the  $Z$-algebra is $c=2(k-1)/(k+2)$.
We have $c=1/2$ for $k=2$ and $c=1$ for $k=4$, for example.
This indicates a possibility that lattice models 
with $c=1/2$ or $c=1$ might be described by the deformed Virasoro
algebra.
Below, we shall discuss that
the representation theories for $k=2$ or $k=4$
do provide us with the physical space of 
the Ising or a special case of the 
eight-vertex model (at $p=q^3$).

Introduce the bosons for the deformed Virasoro algebra as
\begin{eqnarray}
&&[\lambda_n,\lambda_m]=-{1\over n}
{(1-q^{-n})(1-(-)^n q^{{k+2\over 2}n})\over 
1+(-)^n q^{{k\over 2}n}}\delta_{n+m,0},\\
&&\left[\lambda_0,Q\right]=k+2.
\end{eqnarray}
The deformed Virasoro current $T(\zeta)$ is bosonized as
\begin{eqnarray}
&&T(\zeta)=
\Lambda_+((-)^{{1\over 2}} q^{-{k\over 4}}\zeta)
+\Lambda_-((-)^{-{1\over 2}}q^{{k\over 4}}\zeta),\\
&&\Lambda_\pm(\zeta)=:\exp\left(\pm\sum_{n\neq 0}\lambda_n
\zeta^{-n}\right):(-)^{\pm{1\over 2}} q^{\mp {1\over 2}\lambda_0}.
\end{eqnarray}
Two screening currents are defined by
\begin{eqnarray}
S_+(\zeta)\!\!\!\!&=&\!\!\!\!
:\exp\!\left(\!-\sum_{n\neq 0}
{1+(-)^n q^{-{k\over 2}n}\over 1-q^{n}}(-)^{-{n\over 2}}q^{{k+2\over 4}n}
\lambda_n
\zeta^{-n}\right)\!\!: e^Q\zeta^{\lambda_0}\zeta^{k+2\over 2},\\
S_-(\zeta)\!\!\!\!&=&\!\!\!\!
:\exp\!\left(\!+\sum_{n\neq 0}
{1+(-)^n q^{{k\over 2}n}\over 1-(-)^n q^{{k+2\over 2}n}}q^{n\over 2}
\lambda_n
\zeta^{-n}\right)\!\!: e^{-{2\over k+2}Q}\zeta^{-{2\over k+2}\lambda_0}
\zeta^{2\over k+2}.
\end{eqnarray}

Let us consider the Fock spaces 
\begin{eqnarray}
{\cal F}_r={\bf C}[\lambda_{-1},\lambda_{-2},\cdots]
e^{{r\over k+2}Q}|0\rangle\qquad\qquad  (r+{k\over 2}\in {\bf Z}).
\end{eqnarray}
Then, 
the screening charge 
\begin{eqnarray*}
Q=\oint{d\xi\over 2\pi i\xi}S_+(\xi):{\cal F}_r\longrightarrow 
{\cal F}_{r+k+2}
\end{eqnarray*}
is well defined and we have the properties
\begin{eqnarray}
&&[T(\zeta),Q]_+=0,\\
&&QQ=0.
\end{eqnarray}
Thus we obtain the cochain complex
\begin{eqnarray}
\cdots
\mathop{\longrightarrow}^Q {\cal F}_{r-k-2}
\mathop{\longrightarrow}^Q {\cal F}_r
\mathop{\longrightarrow}^Q {\cal F}_{r+k+2}
\mathop{\longrightarrow}^Q 
\cdots.
\end{eqnarray}
The zero-the cohomology $H_r^0$ remains nontrivial and all the other
cohomologies vanish.

Let us study the character of the cohomology
${\rm ch}_{{H}_r^0}={\rm tr}_{{H}_r^0}( x^{-\rho})$.
Here the grading operator is defined by
\begin{eqnarray}
\rho=\sum_{n=1}^\infty   
{n^2(1+(-)^n q^{{k\over 2}n})\over (1-q^{-n})(1-(-)^n q^{{k+2\over 2}n})}
\lambda_{-n}\lambda_n
-{\lambda_0^2-1\over 2(k+2)}-{1\over 8}.
\end{eqnarray}
After taking the alternating sum, we have
\begin{eqnarray}
{\rm ch}_{{H}_r^0}={\rm tr}_{{H}_r^0}( x^{-\rho})
=
x^{{r^2-1\over 2(k+2)}+{1\over 8}}
{\Theta_{x^{k+2}}(x^{r+{k+2\over 2}})\over (x;x)_\infty}.
\end{eqnarray}

The following examples are important for our task:
\begin{eqnarray}
(k=1,c=0)&&{\rm ch}_{{H}_{\pm 1/2}^0}=1,\\
(k=2,c={1\over 2})&&{\rm ch}_{{H}_{\pm 1}^0}=x^{1/8}(-x^2;x^2)_\infty
\qquad(\mbox{Ramond}),\\
     &&{\rm ch}_{{H}_{0}^0}= (-x;x^2)_\infty
\qquad(\mbox{Neveu-Schwarz}),\\
(k=4,c=1)&&{\rm ch}_{{H}_{\pm 1}^0}=x^{1/8}{1\over (x;x^2)_\infty}
\qquad(\mbox{eight-vertex}).
\end{eqnarray}

%%%%%%%%%%%%%%%%%%%%%%%%%%%%%%%%%%%%%%%%%%%%%%%%%%%%%%%%%%%%%
\subsection{bosonic construction at $p=q^3$ from $k=4$ DVA}

When $p=q^{2+2/l}$ $(l=1,2,3,\cdots)$, $R^*(\zeta)$
becomes (anti) diagonal matrix. Then,
the commutation relation for $\Psi(\zeta)$
becomes simple.
Among these, the case $p=q^3$ ($l=2$) is the simplest one, from 
the bosonization point of view.
The case $p=q^4$ ($l=1$) will be treated in the next subsection.
At this moment,
all the other cases ($0\leq p<q^2,p\neq q^3,\neq q^4$) remain unclear.

Let us consider the deformed Virasoro algebra with $k=4$.
We should change notation as $q\rightarrow q^{1/2}$.
After this change, our $q^{1/2}$ for the deformed Virasoro algebra
coincides with  the $q^{1/2}$ for the eight-vertex model.

The type II vertex operator is realized on the space ${H}_1^0$ as
\begin{eqnarray}
&&
\Psi(\zeta)=
((-)^{1\over2}q^{ -{1\over4}}+
       (-)^{-{1\over2}}q^{{1\over4}})^{-1}
T(\zeta).
\end{eqnarray}
Using this and the relation (\ref{phipsi}), we get the formula for $\Phi(\zeta)$.
The type I vertex operator acting on the space $H_1^0$ is realized as
\begin{eqnarray}
&&\Phi(\zeta)
= {(-q^{1\over2};-q^{3\over2})_\infty 
(-q^{3\over2};-q^{3\over2})_\infty\over
 (+q^{1\over2};-q^{3\over2})_\infty
 (+q^{3\over2};-q^{3\over2})_\infty}\times \nonumber\\
&&\qquad\times
\oint{d\xi\over 2\pi i \xi}\phi(\zeta)S_-(\xi)
{\Theta_{-q^{3\over 2}}(-(-)^{-{1\over 2}}q^{1\over 4}\zeta/\xi)\over
\Theta_{-q^{3\over 2}}((-)^{-{1\over 2}}q^{1\over 4}\zeta/\xi) },
\end{eqnarray}
where
\begin{eqnarray}
\phi(\zeta)=
:\exp\left(-\sum_{n\neq 0}
{1+(-)^n q^n\over 1-(-)^n q^{{3\over 2}n}}q^{n\over 4}
\lambda_n
\zeta^{-n}\right): e^{{1\over 3}Q}
\zeta^{{1\over 3}\lambda_0}\zeta^{1\over 3}.
\end{eqnarray}
Thus we arrived at the realization at $p^{1/2}=-q^{3/2}$.
The corner Hamiltonian $D$ is given by
\begin{eqnarray}
D=-\rho-{1\over 8}.
\end{eqnarray}

We can check that 
the operator $t(\zeta)$ satisfying (\ref{comfort})
is obtained both from fusing type II VO's and type I VO's
as (\ref{t}). We omit the detail.

We note that if we further introduce the sign 
change $q^{1/2}\rightarrow -q^{1/2}$,
we get a formula which works for $p^{1/2}=q^{3/2}$.

%%%%%%%%%%%%%%%%%%%%%%%%%%%%%%%%%%%%%%%%%%%%%%%%%%%%%%%%%%%
\subsection{bosonic construction at $p=q^4$ from $k=2$ DVA}

Let us first introduce the vertex operators for 
the Ising model. We closely follow the description given in \cite{FJMMN}.
The particles can be described by 
the Neveu-Schwarz fermions $\psi^{\rm NS}_n$ 
$(n:\mbox{odd integers})$ and the Ramond fermions
$\psi^{\rm R}_n$ $(n:\mbox{even integers})$ satisfying
anti commutation relations
\begin{eqnarray}
&&[\psi^{\rm NS}_n,\psi^{\rm NS}_m]_+=-(q^n+q^{-n})\delta_{n+m,0},\\
&&[\psi^{\rm R}_n,\psi^{\rm R}_m]_+=(q^n+q^{-n})\delta_{n+m,0},
\end{eqnarray}
and write $\psi^{\rm NS}(\zeta)=\sum_{n:{\rm odd}}\psi^{\rm NS}_n
\zeta^{-n}$ and $\psi^{\rm R}(\zeta)=\sum_{n:{\rm even}}\psi^{\rm R}_n
\zeta^{-n}$.
Denote the fermionic vacuums as 
$|{\rm vac}\rangle_{\rm NS}, |{\rm vac}\rangle_{\rm R}$.
We  The type I operators
$\Phi_{\rm NS\pm}^{\rm R}(\zeta)$ and 
$\Phi_{\rm NS\pm}^{\rm R}(\zeta)$ should satisfy
\begin{eqnarray}
&&
\Phi_{\rm NS\pm}^{\rm R}(\zeta_1)\psi^{\rm NS}(\zeta_2)=
\zeta 
{(q^3\zeta^2;q^4)_\infty\over
 (q\zeta^2;q^4)_\infty}
{(q\zeta^2;q^4)_\infty\over
 (q^3\zeta^2;q^4)_\infty}
\psi^{\rm R}(\zeta_2)\Phi_{\rm NS\pm}^{\rm R}(\zeta_1),\\
&&
\Phi_{\rm R}^{\rm NS\pm}(\zeta_1)\psi^{\rm R}(\zeta_2)=
\zeta^2
{(q^3\zeta^2;q^4)_\infty\over
 (q\zeta^2;q^4)_\infty}
{(q\zeta^2;q^4)_\infty\over
 (q^3\zeta^2;q^4)_\infty}
\psi^{\rm NS}(\zeta_2)\Phi_{\rm R}^{\rm NS\pm}(\zeta_1),
\end{eqnarray}
where $\zeta=\zeta_2/\zeta_1$.

Normalize the vertex operators by the conditions
${}_{\rm R}\langle{\rm vac}|
\Phi_{\rm NS+}^{\rm R}(\zeta)
 |{\rm vac}\rangle_{\rm NS}=1$
and
${}_{\rm NS}\langle{\rm vac}|
\Phi_{\rm R}^{\rm NS+}(\zeta)
 |{\rm vac}\rangle_{\rm R}=1$.
We have the parity relations
$\Phi_{\rm NS\pm}^{\rm R}(-\zeta)=\pm \Phi_{\rm NS\pm}^{\rm R}(\zeta)$,
$\Phi_{\rm R}^{\rm NS\pm}(-\zeta)=\pm \Phi_{\rm R}^{\rm NS\pm}(\zeta)$.

The type I vertex operators satisfy the commutation relation
\begin{eqnarray}
&&
w_{\sigma\sigma'}(\zeta_2/\zeta_1)
\Phi_{\rm R}^{\rm NS\sigma}(\zeta_1)
\Phi_{\rm NS\sigma'}^{\rm R}(\zeta_2)
=
\Phi_{\rm R}^{\rm NS\sigma}(\zeta_2)
\Phi_{\rm NS\sigma'}^{\rm R}(\zeta_1),\\
&&
\sum_{\sigma'}
\overline{w}_{\sigma\sigma'}(\zeta_2/\zeta_1)
\Phi_{\rm NS\sigma'}^{\rm R}(\zeta_1)
\Phi_{\rm R}^{\rm NS\sigma'}(\zeta_2)
=
\Phi_{\rm NS\sigma}^{\rm R}(\zeta_2)
\Phi_{\rm R}^{\rm NS\sigma}(\zeta_1),
\end{eqnarray}
where the (normalized) Boltzmann weights for the Ising model
are
\begin{eqnarray}
&&
w_\pm(\zeta)
=
\zeta^{1\mp 1\over 2}
(q^2 \zeta^{\pm 2};q^8)_\infty(q^6 \zeta^{\mp 2};q^8)_\infty\times \\
&&\qquad\qquad\qquad\times
{(q^4\zeta^2;q^4,q^4)_\infty\over (q^2\zeta^2;q^4,q^4)_\infty}
{(q^6\zeta^{-2};q^4,q^4)_\infty\over (q^4\zeta^{-2};q^4,q^4)_\infty},\nonumber\\
&&
\overline{w}_\pm(\zeta)
=
{(q^2;q^4)_\infty\over (q^4;q^8)_\infty^2}
w_\pm(-q/\zeta).
\end{eqnarray}

Setting $k=2$ for the deformed algebra, we have a realization of the Ising model.
First, fermions are realized by $T(\zeta)$ as
\begin{eqnarray}
((-)^{1\over2}q^{ -{1\over2}}+
       (-)^{-{1\over2}}q^{ {1\over2}})^{-1} \;T(\zeta)
=
\left\{
\begin{array}{ll}
\psi^{\rm NS}(\zeta)\qquad &(\mbox{on }{H}_0^0) \\[4mm]
\psi^{\rm R}(\zeta)  &(\mbox{on }{H}_1^0)
\end{array}\right. .
\end{eqnarray}
Then the type I VO's can be obtained 
from the relation (\ref{phipsi}).
If we combine the type I operators as
\begin{eqnarray}
&&
\Phi_{\rm NS}^{\rm R}(\zeta)=\zeta^{1/8}\left(
\Phi_{\rm NS+}^{\rm R}(\zeta)+
\Phi_{\rm NS-}^{\rm R}(\zeta)\right),\\
&&
\Phi_{\rm R}^{\rm NS}(\zeta)=\zeta^{-1/8}\left(
\Phi_{\rm R}^{\rm NS+}(\zeta)+
\Phi_{\rm R}^{\rm NS-}(\zeta)\right),
\end{eqnarray}
we can write down the bosonic expression in a simple manner as
\begin{eqnarray}
&&
\Phi_{\rm NS}^{\rm R}(\zeta)=\phi(\zeta),\\ 
&&
\Phi_{\rm R}^{\rm NS}(\zeta)
= {(-q;-q^2)_\infty (-q^2;-q^2)_\infty\over
 (+q;-q^2)_\infty (+q^2;-q^2)_\infty}\times\\
&&\qquad\times
\oint{d\xi\over 2\pi i \xi}\phi(\zeta)S_-(\xi)
{\Theta_{-q^2}(-(-)^{-{1\over 2}}q^{1\over 2}\zeta/\xi)\over
\Theta_{-q^2}((-)^{-{1\over 2}}q^{1\over 2}\zeta/\xi) },\nonumber
\end{eqnarray}
where
\begin{eqnarray}
\phi(\zeta)=
:\exp\left(-\sum_{n\neq 0}
{1\over 1-(-)^n q^{2n}}(-)^{-{n\over 2}}q^{n}
\lambda_n
\zeta^{-n}\right): e^{{1\over 4}Q}
\zeta^{{1\over 4}\lambda_0}\zeta^{1\over 8}.
\end{eqnarray}
These are the following intertwiners
\begin{eqnarray}
&&\Phi_{\rm NS}^{\rm R}(\zeta):{H}_0^0 \rightarrow {H}_1^0,\\
&&\Phi_{\rm R}^{\rm NS}(\zeta):{H}_1^0 \rightarrow {H}_0^0.
\end{eqnarray}
\medskip

Let us consider the eight-vertex model at $p^{1/2}=-q^2$.
At this point, the eight-vertex model decouples into
two independent Ising models.
Therefore, the vertex operators can be 
written in terms of the Ising VO's. Note that
we have two ways to identify 
the irreducible highest weight modules $V(\Lambda_i)$ 
for the eight-vertex model
\begin{eqnarray}
V(\Lambda_i)&\cong& H_0^0\otimes H_1^0\qquad\qquad {\rm or},\\
&\cong&H_1^0\otimes H_0^0. \nonumber
\end{eqnarray}

First, the type II VO is realized as
\begin{eqnarray}
&&\Psi(\zeta)=((-)^{1\over2}q^{ -{1\over2}}+
       (-)^{-{1\over2}}q^{ {1\over2}})^{-1}
\left(
T(\zeta)\otimes {\rm id}+{\rm id}\otimes T(\zeta)
\right),\\
&&\Psi(\zeta):H_i^0\otimes H_{1-i}^0 \longrightarrow 
H_i^0\otimes H_{1-i}^0\qquad \qquad (i=0,1).
\end{eqnarray}
Next, for the type I VO, we have
\begin{eqnarray}
&&\Phi(\zeta)=
\Phi_{\rm NS}^{\rm R}(\zeta)\otimes
\Phi_{\rm R}^{\rm NS}(\zeta):H_0^0\otimes H_1^0 \longrightarrow 
H_1^0\otimes H_0^0,\\
&&\Phi(\zeta)=
\Phi_{\rm R}^{\rm NS}(\zeta)\otimes
\Phi_{\rm NS}^{\rm R}(\zeta):H_1^0\otimes H_0^0 \longrightarrow 
H_0^0\otimes H_1^0.
\end{eqnarray}
The corner Hamiltonian $D$ is given by
\begin{eqnarray}
D=-\rho\otimes {\rm id}-{\rm id}\otimes \rho-{1\over 8}.
\end{eqnarray}

In \cite{coset}, a fermionic realization of the operator $t(\zeta)$
was presented. While we have here a bosonized version
from fusing type II VO's or type I VO's
as (\ref{t}).
It can be examined that $t(\zeta)$ satisfies (\ref{comfort}).

%%%%%%%%%%%%%%%%%%%%%%%%%%%%%%%%%%%%%%%%%%%%%%%
\subsection{free field construction at $p=q^6$}

Another example of free field realization 
for the eight-vertex model can be investigated
by studying a degeneration of the
integration kernel
$F(\xi_1,\xi_2;p^{1/2},q)$ given by (\ref{Ffor2}). 
Setting $p^{1/2}=q^3$, we
have
\begin{eqnarray}
&&h(\zeta)={(-q\zeta;q^3)_\infty \over(-q^2\zeta;q^3)_\infty } ,\\
&&g(\zeta)={(q^2\zeta;q^3)_\infty \over(q\zeta;q^3)_\infty }
=h(-\zeta)^{-1} ,\\
&&F(\xi_1,\xi_2;q^3,q)
 ={1-\xi^2\over (1+q^{-1}\xi)(1+q\xi)}h(\xi),
\end{eqnarray}
where $\xi=\xi_2/\xi_1$.
Once we realize that the factor 
$(1-\xi^2)/(1+q^{-1}\xi)/(1+q\xi)$ 
in nothing but a two point function for
a $q$-free fermion, we can easily guess the whole structure.

Introduce free bosons and fermions satisfying
\begin{eqnarray}
&&[a_n,a_m]=-{1\over n}{1-q^n\over 1-q^{3n}}(-)^n q^n \delta_{n+m,0},\\
&&[a_0,Q]=1,\\
&&[\psi_n,\psi_m]_+=(-)^n(q^n+q^{-n})\delta_{n+m,0}.
\end{eqnarray}
Define $\phi(\zeta)$ and screening currents $S_\pm(\zeta)$ by
\begin{eqnarray}
&&\phi(\zeta)=
:\exp\left(\sum_{n\neq 0}a_n
\zeta^{-n}\right): e^{ Q}
\zeta^{{1\over 3}a_0}\zeta^{1\over 6},\\
&&S_-(\zeta)=\psi(\zeta)
:\exp\left(-\sum_{n\neq 0}(-)^n a_n
\zeta^{-n}\right): e^{ -Q}
\zeta^{-{1\over 3}a_0}\zeta^{1\over 6},\\
&&S_+(\zeta)=\psi(\zeta)
:\exp\left(\sum_{n\neq 0}{1-q^{3n}\over 1-q^n}q^{-n} a_n
\zeta^{-n}\right): e^{ 3Q}
\zeta^{a_0}\zeta^{3\over 2},
\end{eqnarray}
where $\psi(\zeta)=\sum_n \psi_n \zeta^{-n}$.
We define the Fock spaces as
\begin{eqnarray}
{\cal F}_r={\bf C}[a_{-1},a_{-2},\cdots,\psi_{-1},\psi_{-2},\cdots]
e^{({1\over 2}+3 r)Q}|0\rangle.
\end{eqnarray}

The screening charge
\begin{eqnarray}
Q=\oint {d\xi \over 2 \pi i \xi}S_+(\xi):{\cal F}_r\longrightarrow
{\cal F}_{r+1},
\end{eqnarray}
is well defined, and we can check the nill potency condition 
$QQ=0$. Thus the cochain complex
\begin{eqnarray}
\cdots
\mathop{\longrightarrow}^Q {\cal F}_{-1}
\mathop{\longrightarrow}^Q {\cal F}_0
\mathop{\longrightarrow}^Q {\cal F}_{1}
\mathop{\longrightarrow}^Q 
\cdots,
\end{eqnarray}
and whose cohomology are obtained.
The corner Hamiltonian is chosen to be
\begin{eqnarray}
&&D=-\sum_{n=1}^\infty   
{n^2(1- q^{3n})\over 1-q^{n}}(-)^n q^{-n}
a_{-n} a_n \\
&&
\qquad\qquad+
\sum_{n=1}^\infty   
{n \over q^{n}+q^{-n}}(-)^n 
\psi_{-n} \psi_n
+{4 a_0^2-1\over 24}.\nonumber
\end{eqnarray}

Thus we arrive at the bosonic realization of the type I vertex operator
for the eight-vertex model at $p^{1/2}=q^3$:
\begin{eqnarray}
&&\Phi(\zeta)
= {(+q^{2};q^{3})_\infty 
(+q^{3};q^{3})_\infty\over
 (-q^{2};q^{3})_\infty
 (-q^{3};q^{3})_\infty}\times\\
&&\qquad\times
\oint{d\xi\over 2\pi i \xi}\phi(\zeta)S_-(\xi)
{\Theta_{q^{3}}(- q \zeta/\xi)\over
\Theta_{q^{3}}( q \zeta/\xi) }.\nonumber
\end{eqnarray}

Notice that we have not yet constructed the type II VO at $p=q^6$.
It can be checked that
the operator $t(\zeta)$ satisfying (\ref{comfort})
is obtained by fusing type I VO's
as (\ref{t}).

%%%%%%%%%%%%%%%%%%%%%%%%%%%%%%%%%%
\section{Discussion}
Every free field construction
for $p=q^3,q^4,q^6$ presented in this paper 
has a different structure. 
The number of fields which we need for the construction
differ from each other, and we need to have 
suitable BRST cohomologies for each. 
When we try to generalize our formula,
we eventually have to understand the reason why we have
the basic hypergeometric series in (\ref{Ffor2}) which correspond to
the contraction between screening currents.
At this moment, it is not clear
how such screening operator could be constructed, and
no other examples has been obtained.

To compare our formula (\ref{correlation}) with
Lashkevich and Pugai's one does not seem a straightforward task.
For local operators acting on $n$ adjacent sites,
Lashkevich and Pugai's formula has $n$-fold integral, 
while our formula needs $2n$-fold integral. 
Their formula contains a parameter called
$u_0$ which comes from the intertwining vector, and they had to prove that
the correlation function does not depend on $u_0$.
We do not have such an extra parameter since
we do not use the vertex-face correspondence.
Their formula contains the summations 
arising from the intertwining vectors. 
Our formula does not have these.
Recently, 
Quano found a method to take these summations.
His formula for $p=q^4,q^6$, however, contains more integrals than we have.
We have to have a better understanding of the
correlation functions to go further and study physical quantities
for the eight-vertex model.

\bigskip

\noindent
{\it Acknowledgement.}~~~
The author wish to thank O. Foda, 
Y. Hara, M. Jimbo, N. Kitanine, J.M. Maillet, B.
McCoy,  T. Miwa, S. Odake, Y. Pugai, Y.-H. Quano,  and Y. Yamada
for discussions.

\end{document}